\documentclass[12pt]{article}
\def\R{\hbox{{\rm I}\kern-0.2em{\rm R}\kern0.2em}}
\def\D{\hbox{{\rm I}\kern-0.2em{\rm D}\kern0.2em}}

\def\be{\begin{equation}}
\def\ee{\end{equation}}

\def\({\left(}
\def\){\right)}
\def\[{\left[}
\def\]{\right]}
\def\bc{\begin{center}}
\def\ec{\end{center}}

\begin{document}
{\bf CONSTRUCTING A SPACE FROM THE SYSTEM OF GEODESIC EQUATIONS}

E. FREDERICKS, F. M. MAHOMED, E. MOMONIAT

Centre for Differential Equations, Continuum Mechanics and Applications\\
School of Computational and Applied Mathematics\\
University of the Witwatersrand\\
Wits 2050, South Africa\\E-mail: fmahomed@cam.wits.ac.za

and ASGHAR QADIR

Centre for Advanced Mathematics and Physics\\
National University of Sciences and Technology\\
Campus of the College of Electrical and Mechanical Engineering\\
Peshawar Road, Rawalpindi, Pakistan
and\\
Department of Mathematical Sciences\\
King Fahd University of Petroleum and Minerals\\
Dhahran 31261, Saudi Arabia\\
E-mails: aqadirmath@yahoo.com, aqadirs@comsats.net.pk

\section*{Abstract}

Given a space it is easy to obtain the system of geodesic
equations on it. In this paper the inverse problem of
reconstructing the space from the geodesic equations is addressed.
A procedure is developed for obtaining the metric tensor from the
Christoffel symbols. The procedure is extended for determining if
a second order quadratically semi-linear system can be expressed
as a system of geodesic equations, provided it has terms only
quadratic in the first derivative apart from the second derivative
term. A computer code has been developed for dealing with larger
systems of geodesic equations.

\baselineskip=13.07pt
\section{Introduction}

The development of geometry originates in its applications for map
making [1], but even more from its use in kinematics and dynamics
[2]. As such, it is of interest to look at the interplay of geometry
and dynamics. The path of a test particle in a flat space is a
straight line. The curved space generalization of the straight line
is a geodesic. Thus test particles in curved spaces move along
geodesics. Given a space one easily obtains the system of geodesic
equations on it [3]. In principle it should be easy to obtain the
space from the geodesic equations. This would be of physical
relevance as actual observations would not provide the space but
would provide the paths followed by ``test particles". For example,
in general relativity, one assumes some matter-energy distribution
and then solves the Einstein equations to obtain the metric tensor,
using which one obtains the geodesic equations, which give the paths
for test particles [4]. However, one does not really know the
matter-energy distribution in any actual situation, but only the
observed paths of particles. Consequently, it would be of interest
to be able to determine the metric directly from the geodesic
equations.

Though simple in principle, the problem is complicated by the fact
that the Christoffel symbols are non-linear combinations of the
metric tensor and its first derivative. As such, a system of highly
non-linear first order partial differential equations would have to
be solved to obtain the metric tensor. The problem can be reduced
enormously in complexity by contracting the Christoffel symbols with
the metric tensor, to obtain a system of first order linear partial
differential equations. In general, even this is very complicated to
solve. Further, we would need to check compatibility of the
solutions obtained.

The procedure adopted here uses the skew symmetry of the covariant form of
the Riemann tensor in the first two indices, and the symmetry under
interchange of the first and second pair of indices to provide a system of
linear equations that can be solved simultaneously. If the system does not
decouple we finally have to solve $n$ partial first order differential
equations for one function of $n$ variables. Consequently an arbitrary
constant appears in the solution. If it decouples we need to solve
correspondingly more differential equations and hence more arbitrary
constants appear. These constants generally get determined by inserting the
solutions back into the equations for the metric. In principle it is
possible that they may not be fully evaluated and lead to a class
of metrics. In this paper a specific prescription is provided to reconstruct
the space from the geodesic equations. The metric constructed is unique
(up to some multiplicative constants appearing in the solution). For two
variables the procedure will be explained in detail in the next section.
However, for larger systems the procedure is still too complicated to be
implemented by hand. In section 3 we have provided a general discussion of
the general case and the logic of the computer code to obtain the metric from
the Christoffel symbols. Further, if we only have a system of second order
quadratically semi-linear ordinary differential equations (ODEs) given, that
have only the quadratic term, we would not know whether they could,
consistently, be regarded as a system of geodesic equations. In section 4 we
give a brief discussion of how the code can check whether the system can, or
cannot, be regarded as a system of geodesic equations. In section 5 there are
some specific examples given to illustrate the use of the general procedure. In
the last section we have given a brief summary and discussion of the results.

\section{A System of Two Equations for Two Variables}

The essential principle for obtaining the metric from the Christoffel symbols
may be seen directly by considering a system of two geodesic equations for
two functions of one variable. However, the general procedure involves
additional complications that will be discussed later. For the system of two
equations
\begin{equation}
x'' = a(x,y)x'^2 + 2b(x,y)x'y' + c(x,y)y'^2 ,
\label{eq:it}
\end{equation}
\begin{equation}
y'' = d(x,y)x'^2 + 2e(x,y)x'y' + f(x,y)y'^2 ,
\label{eq:it}
\end{equation}
we can read off the Christoffel symbols as the the negative of the
coefficients of the quadratic terms. Thus
\begin{equation}
\Gamma^1_{11} = -a, \Gamma^1_{12} = -b, \Gamma^1_{22} = -c ,
\Gamma^2_{11} = -d, \Gamma^2_{12} = -e, \Gamma^2_{22} = -f .
\end{equation}
Note that for a general second-order quadratically semi-linear system of
ODEs, the coefficients cannot be assumed to be expressible as Christoffel
symbols. However, if we are given the system of equations as geodesic
equations, we can assume that the coefficients are so expressible. For a
known metric tensor the Christoffel symbols are then given by
\begin{equation}
\Gamma^i_{jk} = \frac{1}{2} g^{im}(g_{jm,k} + g_{km,j} - g_{jk,m}) .
\end{equation}
Now construct the Riemann tensor from these Christoffel symbols
\begin{equation}
R^i_{\;jkl} = \Gamma^i_{jl,k} - \Gamma^i_{jk,l} +\Gamma^i_{ml}
\Gamma^m_{jk} - \Gamma^i_{mk}\Gamma^m_{jl} ,
\end{equation}
where the Einstein summation convention, that repeated indices are summed
over, has been used. Note that the tensor is skew in the last two indices,
$k$ and $l$. As such, when they are equal the tensor is trivially zero. Thus,
without loss of generality, we can set $k = 1, l = 2$. Putting the tensor
into fully covariant form it is skew in its first two indices as well. Using
the metric tensor to lower the index of the curvature tensor, we obtain the
two linear relations for the metric coefficients:
\begin{equation}
g_{11} R^1_{\;112} + g_{12} R^2_{\;112} = 0 ,
\end{equation}
\begin{equation}
g_{12} R^1_{\;212} + g_{22} R^2_{\;212} = 0 .
\end{equation}

There are various possibilities for the Riemann tensor components
being zero or non-zero. Not all possibilities are consistently
allowed. Apart from the case of a flat space, $R^i_{\;jkl} = 0$,
only two possibilities survive: (a) when $i = j$, $R^i_{\;jkl} =
0$; (b) when $i = j$, $R^i_{\;jkl} \not= 0$. In case (b) these
equations can be used to write $g_{11}$ and $g_{22}$ in terms of
$g_{12}$
\begin{equation}
g_{11} = - \frac{R^2_{\;112}}{R^1_{\;112}} g_{12} := A g_{12} ,
\label{case21}
\end{equation}
\begin{equation}
g_{22} = - \frac{R^2_{\;212}}{R^1_{\;212}} g_{12} := B g_{12} .
\label{case22}
\end{equation}

To make the procedure easier to see, write
$g_{11} = p(x,y), g_{12} = q(x,y)$ and $g_{22} = r(x,y)$. In case (b)
using eqs.(7) and (8) we get $p$ and $r$ in terms of $q$. Then, writing the
Christoffel symbols explicitly we obtain the differential equation for $q$
\begin{equation}
q_x = (Ab + a + Bd + e)q , q_y = (Ac + b + Be + f)q .
\end{equation}
The solution for $q$ is provided by integrating eq.(9) relative to $x$ and
$y$ and comparing the arbitrary functions of integration
\begin{equation}
q(x,y) = \alpha (y) exp(\int (Ab + a + Bd + e) dx)
= \beta (x) exp(\int (Ac + b + Be + f) dy) .
\end{equation}
There would appear to be an arbitrary constant still left. This
disappears on using the resulting $p, q, r$ in the expression for the
Christoffel symbols. (Remember that the inverse metric contains the
functions as well as their first derivatives.)

In case (a) $q = 0$. We now get two sets of two partial differential
equations for $p$ and $r$, which can be solved to give
\begin{equation}
p(x,y) = \gamma(y) exp(\int 2a(x,y)dx) = \delta(x) exp(\int 2c(x,y)dy) ;
\end{equation}
\begin{equation}
r(x,y) = \mu (y) exp(\int 2d(x,y)dx) = \nu (x) exp(\int 2f(x,y)dy) .
\end{equation}
There now appear to be two arbitrary constants appearing which are
determined by inserting the expressions back into the Christoffel
symbols. If the constant(s) remain, we obtain a class of metrics for the
same geodesics.

Note the remarkable fact that the algebraic symmetry properties of
the `geometric' intrinsic curvature tensor, from the `ODE' point
of view, are just the compatibility criteria for being able to
obtain the metric from the system of geodesic equations.

 Finally, there remains the case that the space is flat. The
metric tensor can certainly be set as the flat space metric tensor
in Cartesian coordinates and hence the metric is ``reconstructed".
However, this would not be the metric tensor in the coordinates
used. We could now solve the full system of six linear first order
partial differential equations for the three functions $p, q, r$
of two variables $x, y$. The compatibility is now guaranteed.
There are other, neater, methods available as well [5].

\section{The General Procedure}

Whereas, in principle there is nothing new when we have more than
two variables, the problem arises because there are many more
possibilities now. To see how the equations proliferate, consider
the three variable case. We now have $18$ Christoffel symbols for
a system of three geodesic equations. These lead to six
independent components of $R_{ijkl}$. There are now three sets of
constraining equations each of which has three possible index
choices. As such, we have $9$ linearly dependent equations for the
$6$ metric coefficients. If all the components are non-zero, we
have enough equations to be able to obtain the metric coefficients
from three differential equations. In fact if we have five
distinct components non-zero, we could solve the system. However,
we have many possibilities between this case and the flat metric.
For 4 variables there are 10 independent metric coefficients, 20
linearly independent components of $R_{ijkl}$ and 40 Christoffel
symbols for a system of four geodesic equations. This appears to
be a very heavily over-determined system and compatibility checks
would become really long.

The proliferation of equations would have rapidly rendered the
problem intractable were it not for the availability of computer
codes to solve such systems, for many more variables, going through
{\it all} possibilities. We have constructed such a computer code
that enables us to solve the problem in full generality. It is given
at: www.cam.wits.ac.za/inverse.

The logic of the code is as follows. We first differentiate the Christoffel
symbols $\Gamma ^i_{jk}$ relative to all the dependent variables and combine
them to form the curvature tensor of {\it valence} [1, 3], i.e. with one upper
and three lower indices, $R^i_{\;jkl}$. Next we use the symmetry properties of
the covariant form of the Riemann tensor, namely
\be
g_{im}R^m_{\;jkl} = -g_{jm}R^m_{\;ikl}
\label{eb1} ,
\ee
and
\be
g_{im}R^m_{\;jkl} = g_{km}R^m_{\;lij} .
\label{eb2}
\ee
Since (14) is skew in $i, j$, there are $n^3(n - 1)/2$ linearly independent
equations. Further, (15) are $n^4$ equations. There are only $n(n + 1)/2$
independent components of $g_{ij}$. As such, the system must be grossly
over-determined. However, if the $\Gamma ^i_{jk}$ are, indeed, Christoffel
symbols, they must be consistent. As such, one can use the first $n(n+1)/2-1$
of them to obtain all the $g_{ij}$ in terms of one of them (say $g_{11}$). It
is to be noted that since the system is homogeneous, there can be no
non-trivial determination for {\it all} the metric coefficients from here. One
now writes the equation for the metric tensor in terms of the given Christoffel
symbols as
\be
g_{ik,j} + g_{jk,i} - g_{ij,k} = 2 g_{il}\Gamma ^l_{jk} .
\ee
With the full set of equations for all $n^2(n + 1)/2$ independent Christoffel
symbols we can reduce the equations to a system of $n$ first order linear
partial differential equations for one function (say $g_{11}$) of $n$ variables.
We can now solve these and obtain the full metric tensor.

It may happen that the system of equations has a rank less than
$n(n+1)/2-1$. If it is one less, we will need to solve the system
for $n(n+1)/2-2$ components of the metric tensor and then solve
partial differential equations for the last {\it two} metric
coefficients. Similarly, if the rank is $p$ less, we would solve the
system for $n(n+1)/2-p-1$ and then solve the remaining partial
differential equations for those $p$ components. Note that there
would be no need to re-check compatibility of the solutions, other
than to determine the arbitrary constants arising from the solution
of the differential equations.

\section {Consistency Criteria for Systems of
Geodesic Equations}

So far we have taken it for granted that the system given is for
geodesics. It may happen that one obtains a system of equations that
look formally like the system of geodesic equations, in that they
can be written as \be \ddot x^i+\Gamma^i_{jk}\dot x^j\dot x^k=0 ,
\ee but that they cannot be regarded as a system of geodesic
equations. The point is that there is no check that the system of
partial differential equations (16) is internally consistent. To
check this we would require that (14) and (15) form a consistent set
of equations. {\it This is still not enough}! We also need to check
that the first Bianchi identities are satisfied, namely \be
R^i_{\;jkl} + R^i_{\;klj} + R^i_{\;ljk} = 0 . \ee If these are
satisfied we can, indeed, regard the given system of equations as a
system of geodesic equations and possibly use the results of
theorems on global linearizability of the system to obtain the
solution [5].

The computer code we have prepared can be used to obtain the metric tensor if
the system is known to be of geodesic equations and can be used to check the
consistency of the system as geodesic equations.

\section{Computation}

The algorithm is implemented as follows: Specify the order of the
geodesic equation via $\mathbf{n}$.  Assume $g_{ij}=g_{ji}$ and
$R^{i}_{jnn}=0$ by enforcing the following rules
\be\mathbf{SetAttributes[g, Orderless];}\ee\be\mathbf{ R[i\_, j\_,
k\_, l\_] := 0 ; i == j \&\& k == l;}\ee in Mathematica.  We next
introduce the lists {\textbf{SkewSymmetry}}  and {\textbf{Symmetry}}
which caters for all the possible combinations of {\textbf{i,j,k}}
and {\textbf{l}}  when summing over the repeated index in
(\ref{eb1}) and (\ref{eb2}). The independent variables are
represented by {\textbf{X}} in $\mathbf{eq1[i\_, j\_, k\_, l\_]}$
and $\mathbf{eq2[i\_, j\_, k\_, l\_]}$.  The matrices which are then
formed when mapping $\mathbf{eq1[i\_, j\_, k\_, l\_]}$ and
$\mathbf{eq2[i\_, j\_, k\_, l\_]}$ to their corresponding lists
{\textbf{Symmetry}} and {\textbf{SkewSymmetry}} , are stored in
{\textbf{SYM}}  and {\textbf{SkewSYM}}  respectively.  We construct
the $g_{ij}$ metric tensors with {\textbf{gcomponents}} , which is
then used in conjunction with {\textbf{CoefficientArrays}}  to
construct the matrices {\textbf{ASym}}  and {\textbf{ASkewSym}}.  By
choosing $n(n+1)/2-1$ rows from each using the input from
{\textbf{ChooseEqns}} these two matrices are used to form the matrix
{\textbf{A}} and vector {\textbf{b}}.  The vector {\textbf{sol}}
then solves the metric tensors in terms of $g_{11}$ by using
{\textbf{LinearSolve}} in conjunction with the matrix {\textbf{A}}
and vector {\textbf{b}}.  The overdetermined system of linear
partial differential equations, {\textbf{EqnSet16}} , are then used
to solve for $g_{11}$ by  using {\textbf{DSolve}}. The order of the
problem dictates the number of arbitrary functions which will then
have to be solved subsequently.

The implementation of the $n=3$ case is given in the Appendix.

The code has been checked for the following examples.\\
{\bf 1. Systems of two equations:} (a) geodesics on a sphere; (b) a
linearizable system given in [5]; (c) a non-geodesic system (in [5]).\\
{\bf 2. Systems of three equations:} (a) flat space; (b) a 3-sphere;
(c)
linearizable (in [5]); (d) non-geodesic (in [5]).\\
{\bf 3. Physical four dimensional Lorentzian systems:} (a) the
Reissner-Nordstr\"om system (geodesics for a point charged mass, in
which the charge could be taken to be zero to get the Schwarzschild
system); (b) the Kerr system (for a rotating point mass, in which
the rotation could be taken to be zero to reduce to the
Schwarzschild system).

\section {Summary and Discussion}

We have shown explicitly how to construct the metric from the
geodesic equations. In other words, if we knew the geodesics {\it
globally} we could reconstruct the full {\it manifold} with the
metric on it and if we know them locally we can reconstruct the
metric and hence the space, locally. It is remarkable that the
purely geometric entity measuring curvature should provide, when
looked at from the viewpoint of differential equations, the
compatibility conditions for the system to be regarded as describing
geodesics. The significance of this representation is that it
provides a procedure to reduce the system of equations as follows.
We use the conjecture of [6] that for $m$-dimensional sections of
constant curvature it will have an $so(m+1)$ symmetry algebra. We
can use the geometric information to choose the sections of constant
curvature to decouple the system of geodesic equations (using the
procedure of [6]). These $m$ geodesic equations would then be
completely solved and we would have only a system of $n-m$ coupled
equations to be solved. Notice the heavy utilization of purely
geometric considerations.

We solved the problem explicitly for the case of a system of two
geodesic equations analytically. The only problem arose in the
case of the flat space where we anyhow {\it know} the metric in
Cartesian coordinates. If we want to express the metric in the
given coordinates, it may not be so easy. However, the metric
coefficients are directly determinable by solving decoupled first
order partial differential equations for each of the metric
coefficients. There are also more elegant methods available [5].

For larger systems one needs a computer code. Such a code was
developed and its logic is given here. The code is available at:
www.cam.wits.ac.za/inverse. Some examples illustrate the use of the
code. The code further provides a check for the given system of
equations to be consistently regarded as a system of geodesic
equations.

This approach is of importance as it provides a geometric method for
solving systems of ODEs [5]. Further developments may be possible by
embedding the space in higher dimensional spaces in which they
become larger quadratically semi-linear systems. Attempts to convert
cubically semi-linear systems to quadratically semi-linear systems
in a higher dimension are in progress [7]. This procedure is the
inverse of that adopted by Aminova and Aminov [8], in which they use
the geodetic re-parametrization symmetry ($\partial/\partial s$) to
reduce the system by one dimension. Further, it would be possible to
reduce the order of a higher order system by increasing the number
of variables. Thus we could possibly use the same techniques for
higher order ODEs by embedding in correspondingly higher dimensions.
This technique could also be tried to reduce from higher degree
equations to two. It would be of great interest if the approach
could be extended to PDEs as well.

\section*{Appendix}
We implement our code for $\mathbf{n\;=\;3}$. The skew symmetry of
$R^{i}_{jkl}$ implies that without loss of generality, in both
(\ref{eb5}) and (\ref{eb6}) below we have $(k,l) =
\{(1,2),(1,3),(2,3)\}$, while in (\ref{eb6}) $(i,j) =
\{(1,2),(1,3),(2,3)\}$.  The equations (\ref{eb5}) and (\ref{eb6})
are our skew symmetry and symmetry equations repsectively and are
generalized as \be g_{j1} R^{1}_{ikl} + g_{i1} R^{1}_{jkl} + g_{j2}
R^{2}_{ikl} + g_{i2} R^{2}_{jkl} + g_{j3} R^{3}_{ikl} + g_{i3}
R^{3}_{jkl}=0, \label{eb5} \ee \be g_{i1}R^{1}_{jkl} - g_{k1}
R^{1}_{lij} + g_{i2} R^{2}_{jkl} - g_{k2} R^{2}_{lij}+ g_{i3}
R^{3}_{jkl} - g_{k3} R^{3}_{lij}=0. \label{eb6} \ee From (\ref{eb5})
we obtain the nine equations which is generated by $\mathbf{eq2[i\_,
j\_, k\_, l\_]}$ \be g_{11} R^{1}_{112} + g_{12} R^{2}_{112} +
g_{13} R^{3}_{112}=0, \label{eb7}  \ee \be g_{11} R^{1}_{113} +
g_{12} R^{2}_{113} + g_{13} R^{3}_{113}=0, \label{eb8}  \ee \be
g_{11} R^{1}_{123} + g_{12} R^{2}_{123} + g_{13} R^{3}_{123}=0,
\label{eb9}  \ee \be g_{12} R^{1}_{212} + g_{22} R^{2}_{212} +
g_{23} R^{3}_{212}=0, \label{eb10} \ee \be g_{12} R^{1}_{213} +
g_{22} R^{2}_{213} + g_{23} R^{3}_{213}=0, \label{eb11} \ee \be
g_{12} R^{1}_{223} + g_{22} R^{2}_{223} + g_{23} R^{3}_{223}=0,
\label{eb12} \ee \be g_{13} R^{1}_{312} + g_{23} R^{2}_{312} +
g_{33} R^{3}_{312}=0, \label{eb13} \ee \be g_{13} R^{1}_{313} +
g_{23} R^{2}_{313} + g_{33} R^{3}_{313}=0, \label{eb14} \ee \be
g_{13} R^{1}_{323} + g_{23} R^{2}_{323} + g_{33} R^{3}_{323}=0.
\label{eb15} \ee From (\ref{eb6}) we obtain the six equations which
is generated by $\mathbf{eq1[i\_, j\_, k\_, l\_]}$ \be g_{11} \left(
R^{1}_{213} - R^{1}_{312} \right) + g_{12} \left( R^{2}_{213} -
R^{2}_{312} \right) + g_{13} \left( R^{3}_{213} - R^{3}_{312}
\right)=0, \label{eb16} \ee \be g_{11} R^{1}_{223} + g_{12} \left(
R^{2}_{223} - R^{1}_{312} \right) - g_{22} R^{2}_{312} + g_{13}
R^{3}_{223} - g_{23} R^{2}_{312}=0, \label{eb17} \ee \be g_{11}
\left( R^{1}_{312} - R^{1}_{213} \right) + g_{12} \left( R^{2}_{312}
- R^{2}_{213} \right) + g_{13} \left( R^{3}_{312} - R^{3}_{213}
\right)=0, \label{eb18} \ee \be g_{11} R^{1}_{323} + g_{12} \left(
R^{2}_{323} - R^{1}_{313} \right) - g_{22} R^{2}_{313} - g_{23}
R^{3}_{313} + g_{13} R^{3}_{323}=0, \label{eb19} \ee \be g_{11}
R^{1}_{223} - g_{12} \left( R^{1}_{312} - R^{2}_{223} \right) -
g_{22} R^{2}_{312} + g_{13} R^{3}_{223} - g_{23} R^{3}_{312}=0,
\label{eb20} \ee \be - g_{11} R^{1}_{323} + g_{12} \left(
R^{1}_{313} - R^{2}_{323} \right) + g_{22} R^{2}_{313} + g_{23}
R^{3}_{313} - g_{13} R^{3}_{323}=0. \label{eb21} \ee We choose
(\ref{eb10}), (\ref{eb11}), (\ref{eb15}), (\ref{eb16}) and
(\ref{eb17}) to solve for for $g_{12}$, $g_{13}$, $g_{22}$, $g_{23}$
and $g_{33}$.  This choice is done by $\mathbf{\textbf{ChooseEqns} =
\{4, 5, 9, 10, 11\}}$ in our code, which is then used to relate
these metrics to the $g_{11}$ metric.  Here are the relations
$$ \triangle =  R^3_{212} R^3_{223} \left( R^2_{213} \right)^2 -$$
$$ \left( -R^1_{212} \left( R^3_{312} \right)^2 + R^1_{212} R^3_{213} R^3_{312}+ \left( R^2_{312}
   R^3_{212}+R^2_{212} R^3_{213}\right) R^3_{223} \right. $$
$$\left. +R^2_{223} R^3_{212} \left( R^3_{213}-R^3_{312} \right)\right) R^2_{213}  -R^1_{213} R^2_{212}
\left(R^3_{312} \right)^2+ $$
$$R^3_{213} \left( R^2_{212} R^2_{223} R^3_{213}+R^2_{312}  \left( -R^1_{213} R^3_{212}+R^1_{212} R^3_{213}+R^2_{212} R^3_{223}\right)\right) + $$
$$ R^1_{312} \left( R^2_{213} R^3_{212}-R^2_{212} R^3_{213}\right)
   \left(R^3_{213}-R^3_{312}\right)+ $$
\be \left(R^1_{213} R^2_{312} R^3_{212}+\left(R^2_{212}
   \left(R^1_{213}-R^2_{223}\right)-R^1_{212} R^2_{312}\right) R^3_{213} \right) R^3_{312},
\label{eb22} \ee \be g_{12} = g_{11} \left( \left( R^2_{213}
R^3_{212}-R^2_{212} R^3_{213}\right)
   \left( \left(R^1_{312}-R^1_{213}\right) R^3_{223}+R^1_{223}
   \left(R^3_{213}-R^3_{312}\right)\right) \right)/ \triangle,
\label{eb23} \ee
$$ g_{13}= g_{11} \left(\left( R^2_{312} R^3_{212}-R^2_{212} R^3_{312}\right) R_{213}^2+ \right. $$
$$ \left(-\left( R^2_{212} R^2_{223}+R^1_{212} R^2_{312}\right)
R^3_{213}+R^2_{213} \left( R^2_{223} R^3_{212}+R^1_{212}
R^3_{312}\right)+ \right. $$
$$ \left. R_{312}
 \left( R^2_{212} \left( R^3_{213}+R^3_{312}\right)-\left(R^2_{213}+R^2_{312}\right)
 R^3_{212}\right)\right) R^1_{213}+ $$
$$ R_{312}^2 \left( R^2_{213} R^3_{212}-R^2_{212}
R^3_{213}\right)+ $$
$$  R_{223} \left(R^2_{213}-R^2_{312}\right) \left(R^2_{212}
   R^3_{213}-R^2_{213} R^3_{212}\right)+ $$
$$  R_{312} \left(\left(R^2_{212}
   R^2_{223}+R_{212} R^2_{312}\right) R^3_{213}- \right. $$
\be \left. \left. R^2_{213} \left(R^2_{223}
   R^3_{212}+R_{212} R^3_{312}\right)\right)\right)/\triangle,
\label{eb24}  \ee \be g_{22} = g_{11} \left(R^1_{213}
R^3_{212}-R^1_{212} R^3_{213}\right)
\left(\left(R^1_{213}-R^1_{312}\right)
   R^3_{223}+R^1_{223} \left( R^3_{312}-R^3_{213}\right)\right)/\triangle,
   \label{eb25} \ee
\be g_{23} = -g_{11} \left(R^1_{213} R^2_{212}-R^1_{212}
R^2_{213}\right) \left(\left(R^1_{213}-R^1_{312}\right)
R^3_{223}+R^1_{223} \left(
R^3_{312}-R^3_{213}\right)\right)/\triangle,\label{eb25} \ee
$$ g_{33} = g_{11} \left(\left( R^2_{212} R^2_{323} R^3_{223}+R^1_{323} \left( R^2_{212}
R^3_{312}-R^2_{312} R^3_{212}\right)\right) \right. \left(R^1_{213}
\right)^2+   $$
$$  \left(-R^2_{323} \left(R^1_{212}
   R^2_{213} R^3_{223}+R^1_{223} R^2_{212}
   \left(R^3_{213}-R^3_{312}\right)\right)+ \right. $$
$$  R^1_{323} \left(\left(R^2_{212} R^2_{223}+R^1_{212}
R^2_{312}\right) R^3_{213}-R^2_{213} \left(R^2_{223}
R^3_{212}+R^1_{212}
  R^3_{312}\right)\right)+ $$
$$ \left. R_{312} \left(R^1_{323} \left(\left( R^2_{213}+R^2_{312}\right)
   R^3_{212}-R^2_{212} \left( R^3_{213}+R^3_{312}\right)\right)-R^2_{212}
   R^2_{323} R^3_{223}\right)\right) R^1_{213}+ $$
$$ \left(R^1_{312} \right)^2 R^1_{323} \left(R^2_{212}
   R^3_{213}-R^2_{213} R^3_{212}\right)+R^1_{223} \left(R^1_{323}
   \left( R^2_{213}-R^2_{312}\right) \left(R^2_{213} R^3_{212} \right.  \right.$$
$$  \left. \left. -R^2_{212}
R^3_{213}\right)+R^1_{212} R^2_{213} R^2_{323}
  \left(R^3_{213}-R^3_{312}\right)\right)+ $$
$$  R^1_{312} \left(R^1_{212} R^2_{213} R^2_{323}
   R^3_{223}+R^1_{323} \left(R^2_{213} \left(R^2_{223} R^3_{212}+R^1_{212}
   R^3_{312}\right)- \right. \right. $$
\be \left.  \left. \left. \left(R^2_{212} R^2_{223}+R^1_{212}
R^2_{312}\right)
 R^3_{213}\right)\right)\right)/\triangle.
\label{eb26} \ee Substituting into (16) we obtain a system of
differential equations for the $R^{i}_{jkl}$ in terms of the
Christoffel symbols which is represented in our code by
{\textbf{EqnSet16}}. Imposing Christoffel symbols we can determine
the $R^{i}_{jkl}$ from the system of ordinary diferential equations
and hence the metric.   On solving for $g_{11}$ we get
\begin{eqnarray}
g_{11}\left(x_1,x_2,x_3\right)\;=\;&&0\\
g_{11}\left(x_1,x_2,x_3\right)\;=\;&&
   e^{\int_1^{x_1} \Gamma^1_{11}\left({K_1},x_2,x_3\right) \, d{K_1}}
   \int_1^{x_1} e^{-\int_1^{{K_2}} \Gamma^1_{11}\left({K_1},x_2,x_3\right) \, d{K_1}}\nonumber\\
   &&\left(s_2\left({K_2},x_2,x_3\right) \Gamma^2_{11}\left({K_2},x_2,x_3\right)+s_3\left({K_2},x_2,x_3\right) \Gamma^3_{11}\left({K_2},x_2,x_3\right)\right) \,
   d{K_2}\nonumber\\
&&\;+\;e^{\int_1^{x_1} \Gamma^1_{11}\left({K_1},x_2,x_3\right) \,
d{K_1}}
   c_1\left[x_2,x_3\right]
\end{eqnarray}
We can then subsequently solve for the arbitrary function
$c_1\left[x_2,x_3\right]$.

\section*{Acknowledgements}
One of us (AQ) would like to thank the School of Computational and
Applied Mathematics and the Centre for Differential Equations,
Continuum Mechanics and and Applications, of the University of the
Witwatersrand, Johannesburg, South Africa for their hospitality. EM
acknowledges support received from the National Research Foundation
under grant number 2053745.

\section* {References}

[1] See for example Hogben L, {\it Mathematics for the Million} W.W.
Norton and Co. New York 1938.

[2] Barbour J, {\it The Development of Dynamics}, Oxford University Press, 2001.

[3] See any standard book on Differential Geometry, for example:

[4] For example see, Misner CW, Thorne KS and Wheeler JA, {\it
    Gravitation}, W.H. Freeman and Company 1973.

[5] Mahomed FM and Qadir A, ``Linearization of second order
quadratically semi-linear systems of ordinary differential
equations", {\it Nonlinear Dynamics} (to appear).

[6] Feroze T, Mahomed FM and Qadir A, `The Connection Between
Isometries and Symmetries of Geodesic Equations of the Underlying
Spaces', {\it Nonlinear Dynamics} {\bf 45} (2006) 65-74.

[7] Mahomed FM and Qadir A, ``Linearizability criteria for second
order cubically semi-linear systems of ordinary differential
equations", work in progress.

[8] Aminova AV and Aminov N, \textit{Tensor, N. S}, \textbf{62},
(2000) 65.

\end{document}